\documentclass[12pt,a4paper,reqno]{amsart}
\usepackage[a4paper,left=30mm,right=30mm,top=36mm,bottom=36mm]{geometry}

\usepackage{color}
\usepackage{amssymb}
\usepackage{amsmath}
\usepackage{graphicx}
\usepackage{float}
\usepackage[utf8]{inputenc}
\usepackage[english]{babel}
\usepackage{mathtools}
\usepackage{subcaption}
\usepackage{pgfplots}
\usepackage{grffile}
\pgfplotsset{compat=1.18}

\numberwithin{equation}{section}

\usepackage[]{hyperref}
\usepackage[noabbrev, capitalise, nameinlink]{cleveref}
\crefname{equation}{}{}
\crefname{assumption}{Assumption}{Assumptions}
\crefname{rmrk}{Remark}{Remarks}
\crefformat{equation}{\textup{#2(#1)#3}}

\newtheorem{theorem}{Theorem}[section]    
   
\newtheorem{lemma}{Lemma}[section]

\begin{document}

\title[A Cordes framework for stationary FPK equations]{A Cordes framework for stationary Fokker--Planck--Kolmogorov equations}

\author[T. Sprekeler]{Timo Sprekeler}
\address[T. Sprekeler]{Department of Mathematics, Texas A{\&}M University, College Station, TX 77843, USA.}
\email{timo.sprekeler@tamu.edu}

\date{\today}

\begin{abstract} 
We first review the Cordes condition for nondivergence-form differential operators through the lens of Campanato's theory of near operators. We then survey a recently proposed Cordes framework that guarantees the existence and uniqueness of $L^2$ solutions to stationary Fokker--Planck--Kolmogorov equations subject to periodic boundary conditions, and that allows for the construction of a simple finite element method for its numerical approximation. Finally, we propose a Cordes framework for stationary Fokker--Planck--Kolmogorov-type equations subject to a homogeneous Dirichlet boundary condition.
\end{abstract}

\maketitle

\section{Introduction}

In this paper, we consider the stationary Fokker--Planck--Kolmogorov (FPK) equation and a stationary FPK-type equation on $Y:= (0,1)^n$ in the following settings:

\begin{itemize}
\item[] \hspace{-1cm} \textbullet\; Periodic Cordes-type setting \cite{SSZ25}: Seek a (\textit{very weak}) solution $u\in L^2_{\mathrm{per}}(Y)$ to
\begin{align}\label{intro 1}
\begin{split}
&-D^2:(Au) + \nabla\cdot (bu) = 0\quad\text{in }Y,\\ &\qquad u\text{ is $Y$-periodic},\quad \int_Y u = 1,
\end{split}
\end{align}
\hspace{-0.36cm}where $A\in L^{\infty}_{\mathrm{per}}(Y;\mathbb{R}^{n\times n}_{\mathrm{sym}})$ is uniformly elliptic, $b\in L^{\infty}_{\mathrm{per}}(Y;\mathbb{R}^{n})$, and 
\begin{align*}
\exists\, \delta \in \left(\delta_0,1\right]:\quad\frac{\lvert A\rvert^2 + \lvert b \rvert^2}{(\mathrm{tr}(A))^2} \leq \frac{1}{n-1 + \delta}\quad\text{a.e. in }Y
\end{align*} 
\hspace{-0.36cm}with $\delta_0 = (1+4\pi^2)^{-1}$ if $\| b\|_{L^{\infty}(Y)}\neq 0$, and $\delta_0 = 0$ otherwise.

\medskip
\item[] \hspace{-1cm} \textbullet\; Dirichlet Cordes-type setting: Seek a (\textit{very weak}) solution $u\in L^2(Y)$ to
\begin{align}\label{intro 2}
\begin{split}
-D^2:(Au) + \nabla\cdot (bu) &= f\quad\text{in }Y,\\ u &= 0\quad\text{on }\partial Y,
\end{split}
\end{align}
\hspace{-0.36cm}where $A\in L^{\infty}(Y;\mathbb{R}^{n\times n}_{\mathrm{sym}})$ is uniformly elliptic, $b\in L^{\infty}(Y;\mathbb{R}^{n})$, $f\in L^2(Y)$, and 
\begin{align*}
\exists\, \delta \in \left(\hat{\delta}_0,1\right]:\quad\frac{\lvert A\rvert^2 + \lvert b \rvert^2}{(\mathrm{tr}(A))^2} \leq \frac{1}{n-1 + \delta}\quad\text{a.e. in }Y
\end{align*} 
\hspace{-0.36cm}with $\hat{\delta}_0 = (1+\pi^2)^{-1}$ if $\| b\|_{L^{\infty}(Y)}\neq 0$, and $\hat{\delta}_0 = 0$ otherwise.
\end{itemize}

In both settings, we will study existence and uniqueness of an $L^2$ solution as well as its finite element approximation in the absence of any additional regularity assumptions.

\medskip
Our motivation for studying stationary FPK equations stems from periodic homogenization of nondivergence-form equations, where \eqref{intro 1} arises as the problem for the invariant measure; see, e.g., \cite{BLP11}. In the periodic Cordes-type setting stated above with $\|b\|_{L^{\infty}(Y)} = 0$, it was shown in \cite{Spr24} that for any bounded convex domain $\Omega\subset \mathbb{R}^n$, $g\in L^2(\Omega)$, and $\varphi\in H^2(\Omega)$, the nondivergence-form problem
\begin{align*}
-A\left(\frac{\cdot}{\varepsilon}\right):D^2 v^{\varepsilon} &= g\quad\text{in }\Omega,\\
v^{\varepsilon} &= \varphi\quad\text{on }\partial \Omega
\end{align*} 
has a unique solution $v^{\varepsilon}\in H^2(\Omega)$ for any $\varepsilon > 0$, and, as $\varepsilon \searrow 0$, we have that $v^{\varepsilon}$ converges weakly in $H^2(\Omega)$ to the unique solution $v\in H^2(\Omega)$ of the homogenized problem
\begin{align*}
-\overline{A}:D^2 v &= g\quad\text{in }\Omega,\\
v  &= \varphi\quad\text{on }\partial \Omega,
\end{align*}
where the effective diffusion matrix $\overline{A}\in \mathbb{R}^{n\times n}$ is the symmetric positive definite matrix given by
\begin{align*}
\overline{A} := \int_Y Au
\end{align*}
with $u\in L^2_{\mathrm{per}}(Y)$ denoting the solution to \eqref{intro 1}. 

We refer to \cite{CSS20,GST22,GTY20,JZ23,KL16,Spr24,ST21,SWZ25} for recent developments in periodic homogenization, and to \cite{AFL22,AL17,AS14,GST25,GT22,GT23} for recent developments in stochastic homogenization of nondivergence-form equations. From a numerical homogenization perspective, we refer to \cite{CSS20,FGP24,GSS21,KS22,QST24,Spr24,SSZ25} for nondivergence-form problems and to \cite{BLL24,BFP24,HO10,LLM17,LPS18,ZC23} for divergence-form problems with large drift.

While there is a lot of work on the construction of finite element methods for FPK-type equations with smooth coefficients \cite{BWJ96,BS68,KN06,Lan85,Lan91,MB05}, there are only few publications on the case of discontinuous coefficients \cite{LW23,SSZ25,WW20}.

\medskip
This paper is primarily a survey, complemented by some original extensions and a different perspective. We begin with a review of the classical Cordes condition for nondivergence-form differential operators and its connection to Campanato's theory of near operators (see Section \ref{Sec:2}) \cite{Cam94,Cor56,MPS00,SS14}. In Sections \ref{Sec:3} and \ref{Sec:4}, we discuss the existence and uniqueness of a solution $u\in L^2_{\mathrm{per}}(Y)$ to \eqref{intro 1} by reducing the problem to a simple Lax--Milgram problem (see Section \ref{Sec:3}), and its finite element approximation in the periodic Cordes-type setting stated above (see Section \ref{Sec:4}) \cite{SSZ25}. Finally, in Sections \ref{Sec:5} and \ref{Sec:6}, we extend the ideas from the previous sections to study existence and uniqueness of a solution $u\in L^2(Y)$ to \eqref{intro 2} (see Section \ref{Sec:5}) and its finite element approximation in the Dirichlet Cordes-type setting stated above (see Section \ref{Sec:6}).

\section{Cordes condition and Campanato nearness}\label{Sec:2}

Let $\Omega \subset \mathbb{R}^n$ be a bounded convex domain, and let $A = (a_{ij})\in L^{\infty}(\Omega;\mathbb{R}^{n\times n}_{\mathrm{sym}})$ be uniformly elliptic, i.e.,
\begin{align}\label{ell}
\exists\, \lambda,\Lambda > 0:\quad \lambda I_n \leq A \leq \Lambda I_n\quad\text{a.e. in }\Omega.
\end{align}
In the absence of any additional assumptions on $A$, it is well-known that, when $n>2$, the linear nondivergence-form differential operator 
\begin{align}\label{-A:D2}
(-A:D^2) : H^2(\Omega) \cap H^1_0(\Omega) \rightarrow L^2(\Omega),\qquad v\mapsto -A:D^2 v := -\sum_{i,j = 1}^n a_{ij} \partial^2_{ij} v
\end{align}
is not bijective in general; see, e.g., \cite{MPS00}. This is where the Cordes condition comes into play, that is, a condition that limits the scattering of the eigenvalues of $A$ which in turn guarantees bijectivity of $(-A:D^2)$. In this section, we review some known results from the literature through the lens of Campanato's notion of near maps. 

\subsection{The Cordes condition and its geometric interpretation}

The following condition on $A$ was introduced by Cordes in \cite{Cor56} and is commonly referred to as the Cordes condition:
\begin{align}\label{Cor}
\exists\, \delta \in (0,1]:\quad\frac{\lvert A\rvert^2}{(\mathrm{tr}(A))^2} \leq \frac{1}{n-1 + \delta}\quad\text{a.e. in }\Omega,
\end{align}
where $\lvert A\rvert:= \sqrt{A:A}$. Let us denote the eigenvalues of $A$ by $\lambda_1,\ldots,\lambda_n \in L^{\infty}(\Omega)$ and note that $\lambda\leq \lambda_i\leq \Lambda$ a.e. in $\Omega$ for any $i\in\{1,\ldots,n\}$. Recalling that $\lvert A\rvert^2 = \sum_{i=1}^n \lambda_i^2$ and $\mathrm{tr}(A) = \sum_{i=1}^n \lambda_i$, we see that \eqref{Cor} can then be rewritten as
\begin{align*}
\exists\, \delta \in (0,1]:\quad \cos(\theta) = \frac{\sum_{i=1}^n \lambda_i}{\sqrt{\sum_{i=1}^n \lambda_i^2} \sqrt{n}} \geq \sqrt{1-\frac{1-\delta}{n}}\quad\text{a.e. in }\Omega,
\end{align*}
where $\theta(x)$ denotes the angle between the vectors $(\lambda_1(x),\ldots,\lambda_n(x))$ and $(1,\ldots,1)$ in $\mathbb{R}^n$. Hence, the Cordes condition forces the vector $(\lambda_1(x),\ldots,\lambda_n(x))$ of eigenvalues of $A$ to lie inside a cone of sufficiently small angle with axis $(1,\ldots,1)$. When $n = 2$, this cone coincides with the positive quadrant $(0,\infty)^2$, in which case the Cordes condition is implied by uniform ellipticity. Note that the Cordes condition is not implied by uniform ellipticity when $n>2$, and it becomes increasingly restrictive in high dimensions.

\subsection{Bijectivity of \eqref{-A:D2} under Cordes condition via Campanato nearness}

Some algebra quickly reveals the key consequence of the Cordes condition \eqref{Cor}, namely that
\begin{align*}
\left\lvert I_n - \frac{\mathrm{tr}(A)}{\lvert A\rvert^2} A \right\rvert^2 \leq  1-\delta\quad\text{a.e. in }\Omega,
\end{align*}
where $I_n\in \mathbb{R}^{n\times n}$ denotes the identity matrix. In particular, in view of the Miranda--Talenti estimate
\begin{align}\label{MT}
\|D^2 v\|_{L^2(\Omega)}\leq \|\Delta v\|_{L^2(\Omega)}\quad \forall v\in H^2(\Omega) \cap H^1_0(\Omega),
\end{align}
we see that
\begin{align}\label{near1}
\left\| \Delta v - \frac{\mathrm{tr}(A)}{\lvert A\rvert^2} A:D^2 v\right\|_{L^2(\Omega)}\leq  \sqrt{1-\delta}\,\|\Delta v\|_{L^2(\Omega)}\quad \forall v\in H^2(\Omega) \cap H^1_0(\Omega),
\end{align}
i.e., 
\begin{align*}
-\frac{\mathrm{tr}(A)}{\lvert A\rvert^2} A:D^2\qquad\text{is near}\qquad -\Delta
\end{align*}
as maps $H^2(\Omega) \cap H^1_0(\Omega) \rightarrow L^2(\Omega)$ in the sense of Campanato:
\begin{theorem}[nearness of operators \cite{Cam94}]\label{Thm: near}
Let $L_1,L_2:X \rightarrow Y$ be two maps from a nonempty set $X$ to a real Banach space $(Y,\|\cdot\|)$. Suppose that $L_2$ is bijective and that $L_1$ is near $L_2$, that is, there exist constants $\alpha > 0$ and $K\in [0,1)$ such that
\begin{align*}
\| L_2(v_1) - L_2(v_2) - \alpha\left[ L_1(v_1) - L_1(v_2)\right]\| \leq K \|L_2(v_1) - L_2(v_2)\|\quad \forall v_1,v_2\in X.
\end{align*}
Then, $L_1$ is bijective, i.e., for any $f\in Y$ there exists a unique solution $u\in X$ to 
\begin{align*}
L_1(u) = f.
\end{align*}
Further, for this solution $u$ we have the bound
\begin{align*}
\|L_2(u)-L_2(v)\| \leq \frac{\alpha}{1-K} \|f-L_1(v)\|\quad \forall v\in X.
\end{align*}
\end{theorem}

Applying Theorem \ref{Thm: near} to our setting (recall \eqref{near1}), we find that the operator $(-\frac{\mathrm{tr}(A)}{\lvert A\rvert^2} A:D^2):H^2(\Omega) \cap H^1_0(\Omega) \rightarrow L^2(\Omega)$ is bijective and
\begin{align*}
\|\Delta u\|_{L^2(\Omega)} \leq \frac{1}{1-\sqrt{1-\delta}} \left\|\frac{\mathrm{tr}(A)}{\lvert A\rvert^2} A:D^2 u\right\|_{L^2(\Omega)}\quad \forall u\in H^2(\Omega) \cap H^1_0(\Omega).
\end{align*}
Noting that $\frac{\mathrm{tr}(A)}{\lvert A\rvert^2}\in L^{\infty}(\Omega)$ and $\frac{\mathrm{tr}(A)}{\lvert A\rvert^2} > 0$ a.e. in $\Omega$, we obtain the following:

\begin{theorem}[bijectivity of $(-A:D^2)$ under the Cordes condition \cite{Tal65,SS13}]

Let $\Omega \subset \mathbb{R}^n$ be a bounded convex domain, and let $A\in L^{\infty}(\Omega;\mathbb{R}^{n\times n}_{\mathrm{sym}})$ satisfy \eqref{ell} and \eqref{Cor}. Then, the operator 
\begin{align*}
L:H^2(\Omega) \cap H^1_0(\Omega) \rightarrow L^2(\Omega),\qquad v\mapsto Lv:= -A:D^2 v
\end{align*}
is bijective, and we have the bound
\begin{align*}
\|\Delta u\|_{L^2(\Omega)} \leq \frac{1}{1-\sqrt{1-\delta}}\left\|\frac{\mathrm{tr}(A)}{\lvert A\rvert^2} \right\|_{L^\infty(\Omega)} \|L u\|_{L^2(\Omega)}
\end{align*}
for any $u\in H^2(\Omega) \cap H^1_0(\Omega)$.
\end{theorem}

\subsection{Cordes-type condition for operators with lower-order terms}

A natural follow-up question is how the situation changes when lower-order terms are present. To this end, let us consider the differential operator
\begin{align*}
(-A:D^2 - b\cdot \nabla + c):H^2(\Omega) \cap H^1_0(\Omega) \rightarrow L^2(\Omega),
\end{align*}
where $A\in L^{\infty}(\Omega;\mathbb{R}^{n\times n}_{\mathrm{sym}})$ is uniformly elliptic, $b\in L^{\infty}(\Omega;\mathbb{R}^n)$ and $c\in L^{\infty}(\Omega)$ with $c\geq 0$ a.e. in $\Omega$. Consider the following Cordes-type condition from \cite{SS14}:
\begin{align}\label{Cordes-type}
\exists\, (\delta,\lambda) \in (0,1]\times (0,\infty):\quad\frac{\lvert A\rvert^2 + \frac{1}{2\lambda}\lvert b\rvert^2 + \frac{1}{\lambda^2} c^2}{(\mathrm{tr}(A)+\frac{1}{\lambda}c)^2} \leq \frac{1}{n+ \delta}\quad\text{a.e. in }\Omega.
\end{align}
Introducing the renormalization function
\begin{align*}
s := \frac{\mathrm{tr}(A) + \frac{1}{\lambda}c }{\lvert A\rvert^2 + \frac{1}{2\lambda}\lvert b\rvert^2 + \frac{1}{\lambda^2}c^2}\in L^{\infty}(\Omega),
\end{align*}
some algebra quickly reveals the key consequence of the Cordes-type condition \eqref{Cordes-type}, namely that
\begin{align*}
\left\lvert I_n - s A \right\rvert^2 + \frac{1}{2\lambda}\left\lvert  s b \right\rvert^2 + \frac{1}{\lambda^2}\left\lvert \lambda -  s c \right\rvert^2  \leq  1-\delta\quad\text{a.e. in }\Omega.
\end{align*}
In particular, for any $v\in H^2(\Omega) \cap H^1_0(\Omega)$, using that by the Miranda--Talenti estimate \eqref{MT} and integration by parts we have
\begin{align*}
\|D^2 v\|_{L^2(\Omega)}^2 + 2\lambda \|\nabla v\|_{L^2(\Omega)}^2 + \lambda^2 \|v\|_{L^2(\Omega)}^2\leq \|(-\Delta + \lambda) v\|_{L^2(\Omega)}^2,
\end{align*}
we see that
\begin{align*}
\left\| (-\Delta+\lambda) v - ( -sA:D^2 - sb\cdot \nabla + sc)v\right\|_{L^2(\Omega)}\leq  \sqrt{1-\delta}\,\|(-\Delta+\lambda) v\|_{L^2(\Omega)},
\end{align*}
i.e., 
\begin{align*}
-sA:D^2 - sb\cdot \nabla + sc\qquad\text{is near}\qquad -\Delta + \lambda
\end{align*}
as maps $H^2(\Omega) \cap H^1_0(\Omega) \rightarrow L^2(\Omega)$. Noting that the latter map is bijective, we can apply Theorem \ref{Thm: near} to deduce that $-sA:D^2 - sb\cdot \nabla + sc:H^2(\Omega) \cap H^1_0(\Omega) \rightarrow L^2(\Omega)$ is bijective as well, and we have the bound
\begin{align*}
\|(-\Delta + \lambda) u\|_{L^2(\Omega)} \leq \frac{1}{1-\sqrt{1-\delta}}\left\|s( -A:D^2 u - b\cdot \nabla u + c u) \right\|_{L^2(\Omega)}
\end{align*}
for any $u\in H^2(\Omega) \cap H^1_0(\Omega)$. Noting that $s > 0$ a.e. in $\Omega$, we obtain the following:

\begin{theorem}[bijectivity of $(-A:D^2- b\cdot \nabla + c)$ in a Cordes setting \cite{SS14}]

Let $\Omega \subset \mathbb{R}^n$ be a bounded convex domain. Let $A\in L^{\infty}(\Omega;\mathbb{R}^{n\times n}_{\mathrm{sym}})$ satisfy \eqref{ell}, let $b\in L^{\infty}(\Omega;\mathbb{R}^{n})$, and let $c\in L^{\infty}(\Omega)$ with $c\geq 0$ a.e. in $\Omega$. Suppose that the Cordes-type condition \eqref{Cordes-type} holds. Then, the operator 
\begin{align*}
L:H^2(\Omega) \cap H^1_0(\Omega) \rightarrow L^2(\Omega),\qquad v\mapsto Lv:= -A:D^2 v- b\cdot \nabla v + cv
\end{align*}
is bijective, and we have the bound
\begin{align*}
\|(-\Delta + \lambda) u\|_{L^2(\Omega)} \leq \frac{1}{1-\sqrt{1-\delta}}\left\|\frac{\mathrm{tr}(A) + \frac{1}{\lambda}c }{\lvert A\rvert^2 + \frac{1}{2\lambda}\lvert b\rvert^2 + \frac{1}{\lambda^2}c^2} \right\|_{L^\infty(\Omega)} \|L u\|_{L^2(\Omega)}
\end{align*}
for any $u\in H^2(\Omega) \cap H^1_0(\Omega)$.
\end{theorem}

\subsection{Beyond} 

For a more comprehensive overview of elliptic and parabolic equations under Cordes-type conditions, we refer the reader to the books \cite{Lan98,MPS00}. From a numerical perspective, the series of papers \cite{SS13,SS14,SS16} led to extensive subsequent research on the finite element approximation of equations in nondivergence form. 

\section{A Cordes framework for FPK equations: analysis}\label{Sec:3}

We consider the stationary Fokker--Planck--Kolmogorov (FPK) equation on $Y:=(0,1)^n$ subject to periodic boundary conditions and an integral constraint, i.e.,
\begin{align}\label{FPK}
L^{\ast}u := -D^2:(Au) + \nabla\cdot (bu) = 0\quad\text{in }Y,\qquad u\text{ is $Y$-periodic},\qquad \int_Y u = 1,
\end{align}
where $A\in L^{\infty}_{\mathrm{per}}(Y;\mathbb{R}^{n\times n}_{\mathrm{sym}})$ is uniformly elliptic, i.e.,
\begin{align}\label{ell again}
\exists\, \lambda,\Lambda > 0:\quad \lambda I_n \leq A \leq \Lambda I_n\quad\text{a.e. in }Y,
\end{align}
and $b\in L^{\infty}_{\mathrm{per}}(Y;\mathbb{R}^n)$. Note that $L^{\ast}$ is the formal adjoint of the differential operator
\begin{align*}
L := -A:D^2 - b\cdot \nabla.
\end{align*}
We seek solutions to \eqref{FPK} in the space $L^2_{\mathrm{per}}(Y)$. Here, we say that $u\in L^2_{\mathrm{per}}(Y)$ is a (\textit{very weak}) solution to \eqref{FPK} if $\int_Y u = 1$ and
\begin{align*}
(u,L\varphi)_{L^2(Y)} = 0\quad \forall \varphi\in H^2_{\mathrm{per}}(Y).
\end{align*}
In this section, we aim to review a recently proposed Cordes framework for such problems from \cite{SSZ25,Spr24} that guarantees existence and uniqueness of a solution $u\in L^2_{\mathrm{per}}(Y)$ to \eqref{FPK}, and that allows for a simple finite element approximation of $u$ which will be discussed subsequently. 

We set $L^2_{\mathrm{per},0}(Y) := \{ v\in L^2_{\mathrm{per}}(Y): \int_Y v = 0\}$ and $H^k_{\mathrm{per},0}(Y) := H^k_{\mathrm{per}}(Y) \cap L^2_{\mathrm{per},0}(Y)$ for $k\in \mathbb{N}$. 

\subsection{Cordes-type condition}

First, we note that the Cordes-type condition \eqref{Cordes-type} that we previously used to study the operator $(-A:D^2 - b\cdot \nabla + c)$ can never be satisfied when $c = 0$ a.e. due to the fact that $\frac{\lvert M\rvert^2}{(\mathrm{tr}(M))^2}\geq \frac{1}{n}$ for any $M\in \mathbb{R}^{n\times n}$ with $\mathrm{tr}(M)\neq 0$.

Our goal is to impose a condition on the coefficients $A$ and $b$ that guarantees the existence of a function $\gamma\in L^{\infty}_{\mathrm{per}}(Y)$ with $\gamma > 0$ a.e. in $Y$ such that $\gamma L$ is near $-\Delta$ as maps $H^2_{\mathrm{per}}(Y)\rightarrow L^2_{\mathrm{per}}(Y)$. 

This will be achieved by the following Cordes-type condition:
\begin{align}\label{Cor FP}
\exists\, \delta \in \left(\delta_0,1\right]:\quad\frac{\lvert A\rvert^2 + \lvert b \rvert^2}{(\mathrm{tr}(A))^2} \leq \frac{1}{n-1 + \delta}\quad\text{a.e. in }Y
\end{align}
with $\delta_0$ given by
\begin{align}\label{d0}
\delta_0 := (1+4\pi^2)^{-1}\eta,\quad\text{where}\quad \eta:= \begin{cases}
1 &,\text{ if }\|b\|_{L^{\infty}(Y)} \neq 0,\\
0 &,\text{ otherwise.}
\end{cases}
\end{align}
Note that \eqref{Cor FP} reduced to the classical Cordes condition when $\lvert b\rvert = 0$ a.e. in $Y$.

\subsection{Nearness of $L$ to $-\Delta$ after renormalization}\label{Sec:3.2}

Introducing the renormalization function
\begin{align*}
\gamma := \frac{\mathrm{tr}(A)}{\lvert A\rvert^2 + \lvert b\rvert^2} \in L^{\infty}_{\mathrm{per}}(Y),
\end{align*}
some algebra quickly reveals the first key consequence of the Cordes-type condition \eqref{Cor FP}, namely that
\begin{align}\label{id minus gammaA}
\lvert I_n - \gamma A\rvert^2 + \lvert \gamma b\rvert^2 \leq 1-\delta  \quad\text{a.e. in }Y.
\end{align}
In particular, in combination with the Miranda--Talenti-type identity
\begin{align}\label{MT-type}
\|D^2 v\|_{L^2(Y)} = \|\Delta v\|_{L^2(Y)}\quad \forall v\in H^2_{\mathrm{per}}(Y),
\end{align}
and the Poincar\'{e} inequality
\begin{align}\label{Poin}
\|w\|_{L^2(Y)} \leq (2\pi)^{-1} \|Dw\|_{L^2(Y)}\quad\forall w\in H^1_{\mathrm{per},0}(Y;\mathbb{R}^n),
\end{align} 
we see that for any $v\in H^2_{\mathrm{per}}(Y)$ there holds (recall the definition of $\eta$ from \eqref{d0})
\begin{align*}
\|(-\Delta v) - \gamma L v\|_{L^2(Y)} &= \|(\gamma A-I_n):D^2 v + \gamma b\cdot \nabla v\|_{L^2(Y)} \\ &\leq \sqrt{1-\delta} \left(\|\Delta v\|_{L^2(Y)}^2 + \eta  \|\nabla v\|_{L^2(Y)}^2\right)^{1/2}\\
&\leq \sqrt{(1-\delta)\left(1+\eta(2\pi)^{-2}\right)} \;\|\Delta v\|_{L^2(Y)}.
\end{align*}
Therefore, we have the bound
\begin{align}\label{near kappa}
\|(-\Delta v) - \gamma L v\|_{L^2(Y)} \leq \sqrt{1-\kappa} \,\|\Delta v\|_{L^2(Y)}\quad\forall v\in H^2_{\mathrm{per}}(Y),
\end{align}
where
\begin{align}\label{kappa}
\kappa := (\delta - \delta_0)\left(1 + \eta (2\pi)^{-2}\right).
\end{align}
It is at this point that we see the reason we imposed the constraint $\delta\in (\delta_0,1]$ with $\delta_0$ defined in \eqref{d0}. In this way, we guarantee that $\kappa \in (0,1]$.

In view of \eqref{near kappa}, it follows that
\begin{align}\label{near inj}
\gamma L = -\gamma A:D^2 - \gamma b\cdot \nabla \qquad\text{is near}\qquad -\Delta
\end{align}
as maps $H^2_{\mathrm{per}}(Y)\rightarrow L^2_{\mathrm{per}}(Y)$. Let us note that the statement \eqref{near inj} is also true as maps $H^2_{\mathrm{per},0}(Y)\rightarrow L^2_{\mathrm{per}}(Y)$. In particular, as $-\Delta:H^2_{\mathrm{per},0}(Y)\rightarrow L^2_{\mathrm{per}}(Y)$ is injective, Campanato's theory of near operators \cite{Cam94} yields that $\gamma L:H^2_{\mathrm{per},0}(Y)\rightarrow L^2_{\mathrm{per}}(Y)$ is injective, and thus, $L:H^2_{\mathrm{per},0}(Y)\rightarrow L^2_{\mathrm{per}}(Y)$ is injective.

\subsection{Construction of very weak solutions to the FPK problem}

Our goal is to show existence and uniqueness of a solution $u\in L^2_{\mathrm{per}}(Y)$ to \eqref{FPK} in a constructive way that can later serve as a basis for our finite element approximation. Indeed, we will reduce \eqref{FPK} to a Lax--Milgram problem.

\subsubsection{The renormalized FPK problem}

We first show existence and uniqueness of a solution $\tilde{u}\in L^2_{\mathrm{per}}(Y)$ to the renormalized FPK problem
\begin{align}\label{FPK ren}
[\gamma L]^* \tilde{u} := -D^2:(\gamma A\tilde{u}) + \nabla\cdot (\gamma b\tilde{u}) = 0\quad\text{in Y},\qquad \tilde{u}\text{ is $Y$-periodic},\qquad \int_Y \tilde{u} = 1,
\end{align}
i.e., seek $\tilde{u}\in L^2_{\mathrm{per}}(Y)$ such that
\begin{align*}
\int_Y \tilde{u} = 1,\qquad (\tilde{u},\gamma L \varphi)_{L^2(Y)} = 0\quad \forall \varphi\in H^2_{\mathrm{per}}(Y),
\end{align*}
or equivalently,
\begin{align*}
\tilde{u}-1\in  L^2_{\mathrm{per},0}(Y),\qquad (\tilde{u},\gamma L \varphi)_{L^2(Y)} = 0\quad \forall \varphi\in H^2_{\mathrm{per},0}(Y).
\end{align*}
To this end, we make the ansatz
\begin{align*}
\tilde{u} = 1 - \Delta \psi,\qquad \psi\in H^2_{\mathrm{per},0}(Y),
\end{align*}
which is no restriction as $-\Delta:H^2_{\mathrm{per},0}(Y) \rightarrow L^2_{\mathrm{per},0}(Y)$ is bijective. The resulting problem for $\psi$ is the following:
\begin{align}\label{LM}
\text{Seek }\psi\in H^2_{\mathrm{per},0}(Y)\text{ s.t.}\quad (-\Delta \psi,\gamma L\varphi)_{L^2(Y)} = (-1,\gamma L\varphi)_{L^2(Y)}\quad \forall \varphi\in H^2_{\mathrm{per},0}(Y),
\end{align} 
which fits into the framework of the Lax--Milgram theorem as by \eqref{near kappa} we have that
\begin{align*}
(-\Delta v,\gamma Lv)_{L^2(Y)} \geq (1-\sqrt{1-\kappa}) \|\Delta v\|_{L^2(Y)}^2\quad \forall v\in H^2_{\mathrm{per},0}(Y).
\end{align*}
We immediately obtain the following result:
\begin{lemma}[well-posedness and solution structure of the renormalized FPK problem \cite{SSZ25}]\label{Lmm: ren FPK}
Let $Y := (0,1)^n$. Let $A\in L^{\infty}_{\mathrm{per}}(Y;\mathbb{R}^{n\times n}_{\mathrm{sym}})$ and $b\in L^{\infty}_{\mathrm{per}}(Y;\mathbb{R}^n)$ be such that \eqref{ell again} and \eqref{Cor FP} hold, and set $\gamma := \frac{\mathrm{tr}(A)}{\lvert A\rvert^2 + \lvert b\rvert^2}$. Then, there exists a unique solution $\tilde{u}\in L^2_{\mathrm{per}}(Y)$ to the renormalized FPK problem \eqref{FPK ren}. Further, the solution is of the form
\begin{align*}
\tilde{u} = 1-\Delta \psi,
\end{align*}
where $\psi\in H^2_{\mathrm{per},0}(Y)$ is the unique solution to the Lax--Milgram problem \eqref{LM}, and we have that $\tilde{u}\geq 0$ a.e. in $Y$.
\end{lemma}
We omit the proof of nonnegativity of $\tilde{u}$.

\subsubsection{The original FPK problem}

As a direct consequence of Lemma \ref{Lmm: ren FPK} and the positivity of $\gamma$, we obtain the following:
\begin{theorem}[well-posedness and solution structure of the FPK problem \cite{SSZ25}]\label{Thm: FPK}
In the setting of Lemma \ref{Lmm: ren FPK}, there exists a unique solution $u\in L^2_{\mathrm{per}}(Y)$ to the FPK problem \eqref{FPK}. Further, the solution is of the form
\begin{align}\label{r formula}
u = C\gamma \tilde{u} = C\frac{\mathrm{tr}(A)}{\lvert A\rvert^2 + \lvert b\rvert^2} (1-\Delta \psi),
\end{align}
where $\psi\in H^2_{\mathrm{per},0}(Y)$ is the unique solution to the Lax--Milgram problem \eqref{LM}, and $C>0$ is the constant given by
\begin{align*}
C:= (\gamma,\tilde{u})_{L^2(Y)}^{-1} = \left(\frac{\mathrm{tr}(A)}{\lvert A\rvert^2 + \lvert b\rvert^2},1-\Delta \psi\right)_{L^2(Y)}^{-1}.
\end{align*}
\end{theorem}

Note that Theorem \ref{Thm: FPK} is constructive in the sense that it gives a recipe for a numerical scheme to approximate $u$. First, compute $\psi$ using an $H^2_{\mathrm{per},0}(Y)$-conforming finite element method for the Lax--Milgram problem \eqref{LM}. Next, compute the constant $C$ using numerical quadrature. Finally, compute $u$ using the formula given by Theorem \ref{Thm: FPK}. 

\section{A Cordes framework for FPK equations: FE approximation}\label{Sec:4}

Let the setting be as in Theorem \ref{Thm: FPK}, and suppose that $n\in \{2,3\}$. While Theorem \ref{Thm: FPK} immediately suggests a numerical scheme for the approximation of the solution $u\in L^2_{\mathrm{per}}(Y)$ to the FPK problem \eqref{FPK} based on an $H^2$-conforming finite element approximation of an auxiliary function $\psi\in H^2_{\mathrm{per},0}(Y)$, we now aim to construct a numerical scheme that avoids the implementation of an $H^2$-conforming method. 

The key observation is that the formula \eqref{r formula} for $u$ given by Theorem \ref{Thm: FPK} only requires knowledge of $\Delta \psi$ and not of $\psi$ itself. Our goal is to introduce an auxiliary function $\rho\in H^1_{\mathrm{per},0}(Y;\mathbb{R}^n)$ that satisfies $\nabla\cdot \rho = \Delta \psi$ and $\rho$ is the solution to a Lax--Milgram problem in $H^1_{\mathrm{per},0}(Y;\mathbb{R}^n)$. To this end, let us start by recalling the problem for $\psi$ from \eqref{LM}:
\begin{align*}
\text{Seek }\psi\in H^2_{\mathrm{per},0}(Y)\text{ s.t.}\quad (-\Delta \psi,\gamma L\varphi)_{L^2(Y)} = (-1,\gamma L\varphi)_{L^2(Y)}\quad \forall \varphi\in H^2_{\mathrm{per},0}(Y),
\end{align*} 
where $L:=-A:D^2 - b\cdot \nabla$. Let us introduce the differential operator
\begin{align*}
\tilde{L}: H^1_{\mathrm{per},0}(Y;\mathbb{R}^n) \rightarrow L^2_{\mathrm{per}}(Y),\qquad w\mapsto \tilde{L}w := -A:Dw - b\cdot w,
\end{align*}
so that $L\varphi = \tilde{L}(\nabla \varphi)$ for any $\varphi\in H^2_{\mathrm{per}}(Y)$.  

\subsection{Relating the operator $\tilde{L}$ to $-\mathrm{div}$ after renormalization}\label{Sec: 4.1}

We recall from \eqref{id minus gammaA} that the Cordes-type condition \eqref{Cor FP} guarantees that
\begin{align*}
\lvert I_n - \gamma A\rvert^2 + \lvert \gamma b\rvert^2 \leq 1-\delta  \quad\text{a.e. in }Y,\quad\text{where}\quad \gamma = \frac{\mathrm{tr}(A)}{\lvert A\rvert^2 + \lvert b\rvert^2}. 
\end{align*}
In particular, in view of the Poicar\'{e} inequality \eqref{Poin}, we see that
\begin{align}\label{key est mix}
\begin{split}
\| (-\nabla\cdot w) - \gamma \tilde{L} w\|_{L^2(Y)} &= \| (\gamma A - I_n):Dw + \gamma b\cdot w\|_{L^2(Y)} \\ &\leq \sqrt{1-\kappa}\, \|Dw\|_{L^2(Y)}
\end{split}
\end{align}
for any $w\in H^1_{\mathrm{per},0}(Y;\mathbb{R}^n)$, where $\kappa\in (0,1]$ is defined in \eqref{kappa}.

\subsection{The function $\rho$}

In view of \eqref{key est mix}, it is natural to consider the following problem: Seek $\rho\in H^1_{\mathrm{per},0}(Y;\mathbb{R}^n)$ such that
\begin{align}\label{rho}
(-\nabla\cdot \rho,\gamma \tilde{L}w)_{L^2(Y)} + S(\rho,w) = (-1,\gamma \tilde{L}w)_{L^2(Y)}\quad \forall w\in H^1_{\mathrm{per},0}(Y;\mathbb{R}^n),
\end{align}
where the stabilization term $S(\rho,w)$ is chosen as
\begin{align*}
S(\rho,w) := (\mathrm{rot}(\rho),\mathrm{rot}(w))_{L^2(Y)}
\end{align*}
with the convention that $\mathrm{rot}(w):= \partial_2 w_1 - \partial_1 w_2$ when $n = 2$, and $\mathrm{rot}(w):= \nabla \times w$ when $n = 3$. This choice of stabilization is in the spirit of \cite{Gal17b} and motivated by the fact that the first term in \eqref{rho} is close to $(-\nabla\cdot \rho,  -\nabla\cdot w)_{L^2(Y)}$ by \eqref{key est mix}.

We observe that the left-hand side of \eqref{rho} defines a coercive bilinear form on $H^1_{\mathrm{per},0}(Y;\mathbb{R}^n)$. Indeed, in view of \eqref{key est mix} and the identity
\begin{align*}
\|\nabla\cdot w\|_{L^2(Y)}^2 + S(w,w) = \|Dw\|_{L^2(Y)}^2\quad \forall w\in H^1_{\mathrm{per},0}(Y;\mathbb{R}^n),
\end{align*}
we see that
\begin{align*}
(-\nabla\cdot w, \gamma \tilde{L}w)_{L^2(Y)} + S(w,w) \geq \left(1-\sqrt{1-\kappa}\right)\|Dw\|_{L^2(Y)}^2\quad \forall w\in H^1_{\mathrm{per},0}(Y;\mathbb{R}^n).
\end{align*}
Therefore, existence and uniqueness of a function $\rho\in H^1_{\mathrm{per},0}(Y;\mathbb{R}^n)$ solving \eqref{rho} is guaranteed by the Lax--Milgram theorem.

\subsection{Expressing the solution to the FPK problem in terms of $\rho$}

We note that for $\varphi\in H^2_{\mathrm{per}}(Y)$, we can choose $w = \nabla \varphi$ in \eqref{rho} to obtain (recall $\tilde{L}(\nabla \varphi) = L\varphi$)
\begin{align*}
(-\nabla\cdot \rho, \gamma L\varphi)_{L^2(Y)}  =  (-1,\gamma L\varphi)_{L^2(Y)}\quad \forall \varphi\in H^2_{\mathrm{per}}(Y),
\end{align*}
or equivalently,
\begin{align*}
(1-\nabla\cdot \rho, \gamma L\varphi)_{L^2(Y)}  =  0\quad \forall \varphi\in H^2_{\mathrm{per}}(Y),
\end{align*}
i.e., $(1-\nabla\cdot \rho)$ solves the renormalized FPK problem \eqref{FPK ren}. By Lemma \ref{Lmm: ren FPK}, we must have that $\tilde{u} = 1 - \Delta\psi = 1-\nabla\cdot \rho$. In particular, we have shown the following theorem:
\begin{theorem}[another solution formula for the FPK problem \cite{SSZ25}]\label{Thm: second soln form}
Let the setting be as in Theorem \ref{Thm: FPK}, and suppose that $n\in \{2,3\}$. Then, the unique solution $u\in L^2_{\mathrm{per}}(Y)$ to the FPK problem \eqref{FPK} is of the form
\begin{align*}
u = C\gamma \tilde{u} = C\frac{\mathrm{tr}(A)}{\lvert A\rvert^2 + \lvert b\rvert^2} (1-\nabla\cdot \rho),
\end{align*}
where $\rho\in H^1_{\mathrm{per},0}(Y;\mathbb{R}^n)$ is the unique solution to the Lax--Milgram problem \eqref{rho}, and $C>0$ is the constant given by
\begin{align*}
C:= (\gamma,\tilde{u})_{L^2(Y)}^{-1} = \left(\frac{\mathrm{tr}(A)}{\lvert A\rvert^2 + \lvert b\rvert^2},1-\nabla\cdot \rho\right)_{L^2(Y)}^{-1}.
\end{align*}
\end{theorem}

\subsection{Finite element approximation}

The solution formula from Theorem \ref{Thm: second soln form} for the FPK problem \eqref{FPK} suggests a recipe for the numerical approximation of $u$:
\begin{itemize}
\item Step 1: Compute $\rho$ using an $H^1_{\mathrm{per},0}(Y;\mathbb{R}^n)$-conforming finite element method for the Lax--Milgram problem \eqref{rho}. Let us call the resulting approximation $\rho_h$ so that, as $h\searrow 0$, we have $
\|\nabla\cdot (\rho - \rho_h)\|_{L^2(Y)} \longrightarrow 0$.
\item Step 2: Compute $C_h := (\frac{\mathrm{tr}(A)}{\lvert A\rvert^2 + \lvert b\rvert^2},1-\nabla\cdot \rho_h)_{L^2(Y)}^{-1}$, and set
\begin{align*}
u_h := C_h\frac{\mathrm{tr}(A)}{\lvert A\rvert^2 + \lvert b\rvert^2} (1-\nabla\cdot \rho_h).
\end{align*}
Then, we have the error bound $\|u-u_h\|_{L^2(Y)} \lesssim \|\nabla\cdot (\rho - \rho_h)\|_{L^2(Y)}$.
\end{itemize}
 \begin{theorem}[realization of Step 1 \cite{SSZ25}]
Let the situation be as in Theorem \ref{Thm: second soln form}, and let $P_h$ be a closed linear subspace of $H^1_{\mathrm{per},0}(Y;\mathbb{R}^n)$. Then, there exists a unique $\rho_h\in P_h$ such that
\begin{align*}
(-\nabla\cdot \rho_h,\gamma \tilde{L}w_h)_{L^2(Y)} + (\mathrm{rot}(\rho_h),\mathrm{rot}(w_h))_{L^2(Y)} = (-1,\gamma \tilde{L}w_h)_{L^2(Y)}
\end{align*}
for all $w_h\in P_h$, where $\tilde{L}w_h := -A:Dw_h - b\cdot w_h$. Further, we have the error bound
\begin{align*}
\|\nabla\cdot (\rho-\rho_h)\|_{L^2(Y)} \lesssim \inf_{w_h\in P_h} \|D(\rho - w_h)\|_{L^2(Y)}.
\end{align*}
\end{theorem}

\section{Extension to the Dirichlet setting: Analysis}\label{Sec:5}

In this section, our goal is to extend the ideas from Sections \ref{Sec:3} and \ref{Sec:4} to stationary FPK-type equations subject to a homogeneous Dirichlet boundary condition. We consider the problem
\begin{align}\label{Dir prob}
\begin{split}
L^{\ast} u := -D^2:(Au) + \nabla\cdot (bu) &= f\quad\text{in }Y := (0,1)^n,\\
u &= 0\quad\text{on }\partial Y,
\end{split}
\end{align}
where $A\in L^{\infty}(Y;\mathbb{R}^{n\times n}_{\mathrm{sym}})$ satisfies \eqref{ell again} (uniform ellipticity), $b\in L^{\infty}(Y;\mathbb{R}^{n})$, and $f\in L^2(Y)$. We recall that $L^{\ast}$ is the formal adjoint of the differential operator
\begin{align*}
L := -A:D^2 - b\cdot \nabla.
\end{align*} 
We say that $u\in L^2(Y)$ is a (\textit{very weak}) solution to \eqref{Dir prob} if 
\begin{align*}
(u,L\varphi)_{L^2(Y)} = (f,\varphi)_{L^2(Y)}\quad \forall \varphi\in H^2(Y)\cap H^1_0(Y). 
\end{align*} 
In the spirit of Sections \ref{Sec:3} and \ref{Sec:4}, we aim to propose a Cordes framework that guarantees existence and uniqueness of a solution $u\in L^2(Y)$ to \eqref{Dir prob}, and that allows for a simple finite element approximation of $u$ which will be discussed subsequently.

We write $H^1_t(Y;\mathbb{R}^n)$ to denote the subset of $H^1(Y;\mathbb{R}^n)$ consisting of vector fields with vanishing tangential trace on $\partial Y$. 

\subsection{Cordes-type condition}

Our goal is to impose a condition on the coefficients $A$ and $b$ that guarantees the existence of a function $\gamma\in L^{\infty}(Y)$ with $\gamma > 0$ a.e. in $Y$ such that $\gamma L$ is near $-\Delta$ as maps $H^2(Y)\cap H^1_0(Y)\rightarrow L^2(Y)$. 

This will be achieved by the following Cordes-type condition:
\begin{align}\label{Cor FP again}
\exists\, \delta \in \left(\hat{\delta}_0,1\right]:\quad\frac{\lvert A\rvert^2 + \lvert b \rvert^2}{(\mathrm{tr}(A))^2} \leq \frac{1}{n-1 + \delta}\quad\text{a.e. in }Y
\end{align}
with $\hat{\delta}_0$ given by
\begin{align*}
\hat{\delta}_0 := (1+\pi^2)^{-1}\eta,\quad\text{where}\quad \eta:= \begin{cases}
1 &,\text{ if }\|b\|_{L^{\infty}(Y)} \neq 0,\\
0 &,\text{ otherwise.}
\end{cases}
\end{align*}
Note the difference between $\delta_0$ in \eqref{d0} used in the periodic setting, and $\hat{\delta}_0$ used here in the Dirichlet setting.   

\subsection{Nearness of $L$ to $-\Delta$ after renormalization}

Introducing the renormalization function 
\begin{align*}
\gamma := \frac{\mathrm{tr}(A)}{\lvert A\rvert^2 + \lvert b\rvert^2}\in L^{\infty}(Y),
\end{align*}
we have the bound \eqref{id minus gammaA}. Then, arguing as in Section \ref{Sec:3.2}, but with \eqref{MT-type} replaced by the Miranda--Talenti inequality
\begin{align*}
\|D^2 v\|_{L^2(Y)} \leq \|\Delta v\|_{L^2(Y)}\quad \forall v\in H^2(Y)\cap H^1_0(Y),
\end{align*}
and \eqref{Poin} replaced by the Poincar\'{e} inequality
\begin{align}\label{Poin mod}
\| w\|_{L^2(Y)} \leq \pi^{-1} \|Dw\|_{L^2(Y)}\quad \forall w\in H^1_t(Y;\mathbb{R}^n),
\end{align}
we obtain that
\begin{align}\label{near kappah}
\|(-\Delta v) - \gamma L v\|_{L^2(Y)} \leq \sqrt{1-\hat{\kappa}} \,\|\Delta v\|_{L^2(Y)}\quad\forall v\in H^2(Y)\cap H^1_0(Y),
\end{align}
where
\begin{align}\label{kappah}
\hat{\kappa} := (\delta - \hat{\delta}_0)\left(1 + \eta\, \pi^{-2}\right).
\end{align}
Since we assumed $\delta\in (\hat{\delta}_0,1]$, we have that $\hat{\kappa} \in (0,1]$.

In view of \eqref{near kappah}, it follows that
\begin{align*}
\gamma L = -\gamma A:D^2 - \gamma b\cdot \nabla \qquad\text{is near}\qquad -\Delta
\end{align*}
as maps $H^2(Y)\cap H^1_0(Y)\rightarrow L^2(Y)$. In particular, as the latter map is bijective, Theorem \ref{Thm: near} yields that $\gamma L:H^2(Y)\cap H^1_0(Y)\rightarrow L^2(Y)$ is bijective as well, and we have the bound
\begin{align*}
\|\Delta v\|_{L^2(Y)} \leq \frac{1}{1-\sqrt{1-\hat{\kappa}}} \, \|\gamma L v\|_{L^2(Y)} \quad \forall v\in H^2(Y)\cap H^1_0(Y).
\end{align*}
Noting that $\gamma > 0$ a.e. in $Y$, we obtain the following:

\begin{theorem}[bijectivity of $L$]

Let $Y:=(0,1)^n$. Let $A\in L^{\infty}(Y;\mathbb{R}^{n\times n}_{\mathrm{sym}})$ satisfy \eqref{ell again}, let $b\in L^{\infty}(Y;\mathbb{R}^{n})$, and suppose that the Cordes-type condition \eqref{Cor FP again} holds. Then, the operator 
\begin{align*}
L:H^2(Y) \cap H^1_0(Y) \rightarrow L^2(Y),\qquad v\mapsto Lv:= -A:D^2 v- b\cdot \nabla v
\end{align*}
is bijective, and we have the bound
\begin{align*}
\|\Delta v\|_{L^2(Y)} \leq \frac{1}{1-\sqrt{1-\hat{\kappa}}}\left\|\frac{\mathrm{tr}(A)}{\lvert A\rvert^2 + \lvert b\rvert^2} \right\|_{L^\infty(Y)} \|L v\|_{L^2(Y)}
\end{align*}
for any $v\in H^2(Y) \cap H^1_0(Y)$, where $\hat{\kappa}\in (0,1]$ is given by \eqref{kappah}.
\end{theorem}

\subsection{Construction of very weak solutions to the FPK problem}

Our goal is to show existence and uniqueness of a solution $u\in L^2(Y)$ to \eqref{Dir prob} in a constructive way that can later serve as a basis for our finite element approximation. Indeed, we will reduce \eqref{Dir prob} to a Lax--Milgram problem.

\subsubsection{The renormalized problem}

We first show existence and uniqueness of a solution $\tilde{u}\in L^2(Y)$ to the renormalized FPK-type problem
\begin{align}\label{Dir ren}
\begin{split}
[\gamma L]^{\ast} \tilde{u} := -D^2:(\gamma A\tilde{u}) + \nabla\cdot (\gamma b\tilde{u}) &= f\quad\text{in }Y,\\
\tilde{u} &= 0\quad\text{on }\partial Y,
\end{split}
\end{align}
i.e., seek $\tilde{u}\in L^2(Y)$ such that
\begin{align*}
(\tilde{u},\gamma L\varphi)_{L^2(Y)} = (f,\varphi)_{L^2(Y)}\quad \forall \varphi\in H^2(Y)\cap H^1_0(Y). 
\end{align*} 
To this end, we make the ansatz
\begin{align*}
\tilde{u} = - \Delta \psi,\qquad \psi\in H^2(Y)\cap H^1_0(Y),
\end{align*}
which is no restriction as $-\Delta:H^2(Y)\cap H^1_0(Y) \rightarrow L^2(Y)$ is bijective. The resulting problem for $\psi$ is the following:
\begin{align}\label{LM ag}
\text{Seek }\psi\in H:=H^2(Y)\cap H^1_0(Y)\text{ s.t.}\quad (-\Delta \psi,\gamma L\varphi)_{L^2(Y)} = (f, \varphi)_{L^2(Y)}\quad \forall \varphi\in H,
\end{align} 
which fits into the framework of the Lax--Milgram theorem as by \eqref{near kappah} we have that
\begin{align*}
(-\Delta v,\gamma Lv)_{L^2(Y)} \geq (1-\sqrt{1-\hat{\kappa}}) \|\Delta v\|_{L^2(Y)}^2\quad \forall v\in H^2(Y)\cap H^1_0(Y).
\end{align*}
We immediately obtain the following result:
\begin{lemma}[well-posedness and solution structure of the renormalized problem]\label{Lmm: ren Dir}
Let $Y := (0,1)^n$. Let $A\in L^{\infty}(Y;\mathbb{R}^{n\times n}_{\mathrm{sym}})$ and $b\in L^{\infty}(Y;\mathbb{R}^n)$ be such that \eqref{ell again} and \eqref{Cor FP again} hold, and set $\gamma := \frac{\mathrm{tr}(A)}{\lvert A\rvert^2 + \lvert b\rvert^2}$. Then, there exists a unique solution $\tilde{u}\in L^2(Y)$ to the renormalized problem \eqref{Dir ren}. Further, the solution is of the form
\begin{align*}
\tilde{u} = -\Delta \psi,
\end{align*}
where $\psi\in H^2(Y)\cap H^1_0(Y)$ is the unique solution to the Lax--Milgram problem \eqref{LM ag}.
\end{lemma}

\subsubsection{The original problem}

As a direct consequence of Lemma \ref{Lmm: ren Dir} and the positivity of $\gamma$, we obtain the following:
\begin{theorem}[well-posedness and solution structure of the original problem]\label{Thm: Dir}
In the setting of Lemma \ref{Lmm: ren Dir}, there exists a unique solution $u\in L^2(Y)$ to the problem \eqref{Dir prob}. Further, the solution is of the form
\begin{align}\label{u formula}
u = \gamma \tilde{u} = \frac{\mathrm{tr}(A)}{\lvert A\rvert^2 + \lvert b\rvert^2} (-\Delta \psi),
\end{align}
where $\psi\in H^2(Y)\cap H^1_0(Y)$ is the unique solution to the Lax--Milgram problem \eqref{LM ag}.
\end{theorem}

\section{Extension to the Dirichlet setting: FE approximation}\label{Sec:6}

Let the setting be as in Theorem \ref{Thm: Dir}, and suppose that $n\in \{2,3\}$. We now discuss the finite element approximation of the unique solution $u\in L^2(Y)$ to \eqref{Dir prob}. We assume that we are given $F\in H^1(Y;\mathbb{R}^n)$ such that
\begin{align}\label{f struc}
f = -\nabla\cdot F\quad\text{a.e. in }Y.
\end{align}
The existence of $F$ is guaranteed by the surjectivity of $\Delta: H^2(Y)\rightarrow L^2(Y)$.

We begin by observing that the formula \eqref{u formula} for the solution $u$ to \eqref{Dir prob} only requires knowledge of $\Delta \psi$ and not of $\psi$ itself. In the spirit of Section \ref{Sec:4}, we will introduce an auxiliary function $\rho\in H^1_t(Y;\mathbb{R}^n)$ that satisfies $\nabla\cdot \rho = \Delta \psi$ and $\rho$ is the solution to a Lax--Milgram problem in $H^1_t(Y;\mathbb{R}^n)$. To this end, let us start by recalling the problem for $\psi$ from \eqref{LM ag} (and using \eqref{f struc}):
\begin{align*}
\text{Seek }\psi\in H:=H^2(Y)\cap H^1_0(Y)\text{ s.t.}\quad (-\Delta \psi,\gamma L\varphi)_{L^2(Y)} = (F, \nabla\varphi)_{L^2(Y)}\quad \forall \varphi\in H,
\end{align*} 
where $L:= -A:D^2  - b\cdot \nabla$. We introduce the differential operator
\begin{align*}
\tilde{L}:H^1_t(Y;\mathbb{R}^n) \rightarrow L^2(Y),\qquad w\mapsto \tilde{L}w := -A:Dw - b\cdot w,
\end{align*}
so that $L\varphi = \tilde{L}(\nabla \varphi)$ for any $\varphi\in H^2(Y)\cap H^1_0(Y)$.

\subsection{Relating $\tilde{L}$ to $-\mathrm{div}$ after renormalization}

Analogously to Section \ref{Sec: 4.1}, using \eqref{id minus gammaA} and \eqref{Poin mod}, we find that
\begin{align}\label{key est mix mod}
\| (-\nabla\cdot w) - \gamma \tilde{L} w\|_{L^2(Y)} \leq \sqrt{1-\hat{\kappa}}\, \|Dw\|_{L^2(Y)}\quad \forall w\in H^1_t(Y;\mathbb{R}^n),
\end{align}
where $\hat{\kappa}\in (0,1]$ is defined in \eqref{kappah}.

\subsection{The function $\rho$}

We consider the following problem: Seek $\rho\in H^1_t(Y;\mathbb{R}^n)$ such that
\begin{align}\label{rho Dir}
(-\nabla\cdot \rho,\gamma \tilde{L}w)_{L^2(Y)} + (\mathrm{rot}(\rho),\mathrm{rot}(w))_{L^2(Y)} = (F,w)_{L^2(Y)}\quad \forall w\in H^1_{t}(Y;\mathbb{R}^n).
\end{align} 

We observe that the left-hand side of \eqref{rho Dir} defines a coercive bilinear form on $H^1_{t}(Y;\mathbb{R}^n)$. Indeed, in view of \eqref{key est mix mod} and the inequalities
\begin{align*}
\|\nabla\cdot w\|_{L^2(Y)}^2 \leq \|Dw\|_{L^2(Y)}^2 \leq  \|\nabla\cdot w\|_{L^2(Y)}^2 + \|\mathrm{rot}(w)\|_{L^2(Y)}^2 \quad \forall w\in H^1_{t}(Y;\mathbb{R}^n),
\end{align*}
we see that
\begin{align*}
(-\nabla\cdot w, \gamma \tilde{L}w)_{L^2(Y)} + \|\mathrm{rot}(w)\|_{L^2(Y)}^2 \geq \left(1-\sqrt{1-\hat{\kappa}}\right)\|Dw\|_{L^2(Y)}^2\quad \forall w\in H^1_{t}(Y;\mathbb{R}^n).
\end{align*}
Therefore, existence and uniqueness of a function $\rho\in H^1_{t}(Y;\mathbb{R}^n)$ solving \eqref{rho Dir} is guaranteed by the Lax--Milgram theorem.

\subsection{Expressing the solution to \eqref{Dir prob} in terms of $\rho$}

We note that for $\varphi\in H^2(Y)\cap H^1_0(Y)$, we can choose $w = \nabla \varphi$ in \eqref{rho Dir} to obtain (recall $\tilde{L}(\nabla \varphi) = L\varphi$)
\begin{align*}
(-\nabla\cdot \rho, \gamma L\varphi)_{L^2(Y)}  =  (F,\nabla\varphi)_{L^2(Y)} = (f,\varphi)_{L^2(Y)}\quad \forall \varphi\in H^2(Y)\cap H^1_0(Y),
\end{align*}
i.e., $-\nabla\cdot \rho \in L^2(Y)$ solves the renormalized problem \eqref{Dir ren}. By Lemma \ref{Lmm: ren Dir}, we must have that $\tilde{u} = - \Delta\psi = -\nabla\cdot \rho$. In particular, we have shown the following theorem:
\begin{theorem}[another solution formula for \eqref{Dir prob}]\label{Thm: second soln form Dir}
Let the setting be as in Theorem \ref{Thm: Dir}, and suppose that $n\in \{2,3\}$. Then, the unique solution $u\in L^2(Y)$ to the problem \eqref{Dir prob} is of the form
\begin{align*}
u = \gamma \tilde{u} = \frac{\mathrm{tr}(A)}{\lvert A\rvert^2 + \lvert b\rvert^2} (-\nabla\cdot \rho),
\end{align*}
where $\rho\in H^1_{t}(Y;\mathbb{R}^n)$ is the unique solution to the Lax--Milgram problem \eqref{rho Dir}.
\end{theorem}

\subsection{Finite element approximation}

The solution formula from \eqref{Thm: second soln form Dir} for the problem \eqref{Dir prob} suggests a recipe for the numerical approximation of $u$:
\begin{itemize}
\item Step 1: Compute $\rho$ using an $H^1_{t}(Y;\mathbb{R}^n)$-conforming finite element method for the Lax--Milgram problem \eqref{rho Dir}. Let us call the resulting approximation $\rho_h$ so that, as $h\searrow 0$, we have $
\|\nabla\cdot (\rho - \rho_h)\|_{L^2(Y)} \longrightarrow 0$.
\item Step 2: Set
\begin{align*}
u_h = \frac{\mathrm{tr}(A)}{\lvert A\rvert^2 + \lvert b\rvert^2} (-\nabla\cdot \rho_h).
\end{align*}
Then, we have the error bound $\|u-u_h\|_{L^2(Y)} \lesssim \|\nabla\cdot (\rho - \rho_h)\|_{L^2(Y)}$.
\end{itemize}

\begin{theorem}[realization of Step 1]
Let the situation be as in Theorem \ref{Thm: second soln form Dir}, and let $P_h$ be a closed linear subspace of $H^1_{t}(Y;\mathbb{R}^n)$. Then, there exists a unique $\rho_h\in P_h$ such that
\begin{align*}
(-\nabla\cdot \rho_h,\gamma \tilde{L}w_h)_{L^2(Y)} + (\mathrm{rot}(\rho_h),\mathrm{rot}(w_h))_{L^2(Y)} = (F,w_h)_{L^2(Y)}
\end{align*}
for all $w_h\in P_h$, where $\tilde{L}w_h := -A:Dw_h - b\cdot w_h$. Further, we have the error bound
\begin{align*}
\|\nabla\cdot (\rho-\rho_h)\|_{L^2(Y)} \lesssim \inf_{w_h\in P_h} \|D(\rho - w_h)\|_{L^2(Y)}.
\end{align*}
\end{theorem} 

\bibliographystyle{plain}
\bibliography{ref_Spr}

\end{document}